\documentclass{amsart}
 \newtheorem{thm}{Theorem}[section]
 
 \newtheorem{lem}[thm]{Lemma}
 \newtheorem{prop}[thm]{Proposition}
 
 \theoremstyle{definition}
 \newtheorem{defn}[thm]{Definition}
 \theoremstyle{remark}
 \newtheorem{rem}[thm]{Remark}
 \numberwithin{equation}{section}
 \newcommand{\eps}{\varepsilon}
 \newcommand{\To}{\longrightarrow}
 \newcommand{\Real}{{\mathbb{R}}}
 
 \newcommand{\Rd}{{\Real^d}}
 
 \newcommand{\E}{{\bf E}}
 \newcommand{\lap}{{\Delta^{\alpha/2}}}
 \newcommand{\lapN}{{\Delta_N^{\alpha/2}}}
 
 \DeclareMathOperator{\supp}{supp}
 \DeclareMathOperator{\Cap}{Cap}
 \DeclareMathOperator{\diam}{diam}

\begin{document}

 \title[Scattering length for stable processes]{Scattering length 
for stable processes}

 \author{Bart{\l}omiej Siudeja}

\address{Department of Mathematics, Purdue University, West Lafayette, Indiana
47906}

\email{siudeja@math.purdue.edu}

\thanks{Supported in part by NSF Grant \# 9700585-DMS}

\subjclass[2000]{Primary 60G52; Secondary 31C15}

\keywords{stable-like processes, reflected stable process, Neumann eigenvalue}

\begin{abstract}
  Let $\alpha\in(0,2)$ and $X_t$ be a symmetric $\alpha$-stable process. We
  define scattering length $\Gamma(v)$ of a positive potential $v$ and  prove several 
of its basic properties. We use the scattering length to find estimates
  for the first eigenvalue of the Schr\"odinger operator of the ``Neumann'' 
  fractional Laplacian in a cube with potential $v$. 
\end{abstract}

\maketitle

 \section{Introduction}
 The purpose of this paper is to define scattering length for 
 the symmetric $\alpha$-stable processes, or equivalently, for a fractional Laplacian, 
 and to prove some of its basic properties.  
 The scattering length has been studied for Brownian motion and the classical 
 Laplacian by many authors, see \cite{K}, \cite{KL}, \cite{T1}, \cite{T2}. The last two papers
 contain applications. It is possible for example  to give a bound for
the first eigenvalue of the Schr\"odinger operator of the Neumann Laplacian in a cube using
this quantity.
 Scattering length is also important in mathematical physics
 where it arises in many situations,  including the  study of neutron scattering. 

 This paper is the first attempt to define and study scattering length for processes different than Brownian motion. 
 As an application of stable scattering length we prove estimates for the first eigenvalue of the
Schr\"odinger operator of the ``Neumann'' fractional Laplacian in a cube.  This result is similar to the
one obtained in \cite{T1} for the Laplacian. 

 For simplicity we assume that $d>\alpha$, where $d$ is the dimension and $\alpha\in(0,2)$ is the parameter of the process. 
 Throughout the paper all the constants depend on both $d$ and $\alpha$. 
 Dependence on any other parameter will be indicated explicitly. 
 We also adopt the convention that constants may change their values
 from line to line as long as they stay positive.

 We define a symmetric $\alpha$-stable
 process $X_t$ as a Markov process with independent and homogeneous increments
 and characteristic function 
 \begin{equation}
   \E^0(\exp(i\xi X_t))=\exp(-t|\xi|^\alpha).
 \end{equation}
 It is well known that this process has the  generator $\lap = -(-\Delta)^{\alpha/2}$, 
 where $\Delta$ is a classical Laplacian on $\Rd$. For an overview of recent results
 for the potential theory of this process we refer the reader to \cite{C}. 
The  quadratic form for this process 
 is given by
 \begin{equation}
   {\mathcal E}_X(u,u)=\int_\Rd\int_\Rd {(u(x)-u(y))^2\over|x-y|^{d+\alpha}}dxdy.
 \end{equation}
and its domain is  $W^{\alpha/2,2}(\Rd)$, the fractional Sobolev space. See \cite{BBC} for 
details about quadratic forms and domains for the generators of stable processes.

 We can also define the ``Neumann'' fractional Laplacian $\lapN$
 on an open set $\Omega$ as the operator with the quadratic form
 \begin{equation}
   {\mathcal E}_Y(u,u)=\int_\Omega\int_\Omega {(u(x)-u(y))^2\over|x-y|^{d+\alpha}}dxdy
 \end{equation}
 and the domain $W^{\alpha/2,2}(\Omega)$. Here we also refer the reader to \cite{BBC} for details about 
the definition of this operator and for properties of the process
 $Y_t$ generated by it.  In the rest of this paper we shall refer to this 
 process as a reflected stable process in $\Omega$.

To state the main result of this paper we first need a definition of the scattering length.
 Here we just give a quick summary. A precise definition will be given in the section two below.
 Let $v$ be a positive function (assume for now that it is bounded with compact support) and 
 let $U_v$ be the capacitory potential of the function $v$. That is, 
 \begin{equation}
 U_v(x)=\E^x\exp\left(-\int_0^\infty v(X_s)ds\right),
 \end{equation}
 and $\mu_v$ its capacitory measure
 \begin{equation}
 \mu_v=-\lap U_v.
 \end{equation}
 We define the scattering length $\Gamma(v)$ by
 \begin{equation}
 \Gamma(v)=\int_\Rd d\mu_v(x)
 \end{equation}
The main result of the paper is the following
 \begin{thm}\label{main}
   Let $\Omega$ be a cube in $\Rd$, $0\leq v\in L^1(\Omega)$ and 
   $\lambda_1(v)$ be the first eigenvalue of the operator $-\lapN+v$ in $\Omega$ 
   (the Schr\"odinger operator of the ``Neumann'' fractional Laplacian on $\Omega$). 
   Then there exists a constant $C_1(\Omega)$ such that 
   \begin{equation}
     C_1(\Omega)\Gamma(v)\leq \lambda_1(v).
   \end{equation}
   Furthermore,  there exists a constant  $\beta=\beta(\Omega)>0$ such that whenever $\Gamma(v)\leq \beta$, then
   \begin{equation}
     \lambda_1(v)\leq C_2(\Omega)\Gamma(v).
   \end{equation}  
 \end{thm}
 \begin{rem}
   The upper bound is valid for any bounded domain $\Omega$.
 \end{rem}
 
 The idea of the proof is roughly the following.  By choosing an appropriate representative of stable-like
 processes (see [2] and the definition in Section 5) we are able to relate our result to a similar one for Brownian motion
 and for reflected Brownian motion. This allows  us to prove our main theorem for this particular process, and hence
 for any stable-like process, and  in particular for the reflected stable process as in our theorem. 

 The rest of the paper is organized as follows.  In Section 2  we give the precise definition of
scattering length. In Section 3 we prove some properties of capacitory potential and
scattering length. The proofs in this section are easy and they carry over from the Brownian case to the stable
case with minimal changes. We present them here for the sake of completeness. Section 5
contains the necessary estimates required in the proof of  Theorem \ref{main}.  
Sections 4 and 6 contain the proof of Theorem \ref{main}.

 \section{Definitions}
 We will give the probabilistic definition of  scattering length but first 
 we need to define a potential operator for the symmetric stable processes. 
 \begin{defn} For any nonnegative function $f$ define 
   $U$  by
   \begin{gather}
     U[f](x)=\E^x\left(\int_0^\infty f(X_s) ds\right). 
   \end{gather}
   
 \end{defn}
 It is well known that for the symmetric stable processes this operator is given
 by a Riesz kernel (see e.g. \cite{BG})
 \begin{gather}
   U[f](x)=\int_\Rd c{f(y)\over |x-y|^{d-\alpha}}dy.
 \end{gather}
 We have
 \begin{lem}\label{finite}
  If $f\in L^1(\Rd)$, then $U[f](x)$ is finite for almost all $x\in\Rd$. 
  If $f$ is in $L^\infty$, then $U[f](x)$ is finite everywhere.
 \end{lem}
 \begin{proof}
   Set  $g(x)=c_{\alpha,d}/|x|^{d-\alpha}.$
   \begin{gather}
     \begin{split}
       U[f](x)&=(f\circ g)(x)=(f\circ (g 1_{B(0,1)}))(x)+(f\circ (g 1_{B^c(0,1)}))(x).
     \end{split}
   \end{gather}
   The first term is the  convolution of two $L^1$ functions, hence in  $L^1$.
   Since $g\leq 1$ outside the ball
   $B(0,1)$, the  second term is bounded above by $\|f\|_1$.
   Hence whenever $f\in L^1(\Rd)$, the first
   term is bounded by $\|f\|_\infty\|g 1_{B(0,1)}\|_1$, therefore $U[f](x)$ is finite
   everywhere.
 \end{proof}
 Let $v\in L^1(\Rd)$ be positive. Let
 \begin{equation}\label{def1}
   1-U_v^t(x)=e^{t(\lap-v)}1(x)=\E^x\left(e^{-\int_0^t v(X_s)ds}\right).
 \end{equation}
 Using $U_v^t$ we can define the capacitory potential of $v$ by 
 \begin{equation}
   U_v(x)=\lim_{t\To\infty}U_v^t(x)=\E^x\left( 1-e^{-\int_0^\infty v(X_s)ds}\right). 
 \end{equation}

Before we define a scattering length, let us relate it to the capacity of sets.
Let $K$ be a Kac regular set (see \cite{S} for details). Informally, the set $K$ is 
Kac regular if after entering the set $K$ the process will stay there for a positive amount of time.
Put
\begin{gather}
  v_K=\infty\mbox{ on }K\mbox{, and }0\mbox{ outside}.
\end{gather}
Under the assumption of Kac regularity the notion of 
capacitory potential $U_{v_K}$ of $v_K$  coincides with the capacitory 
potential $U[\mu_K]$ of the set $K$, where $\mu_K$ is an equilibrium measure 
on $K$ (see \cite{S})
\begin{gather}
  U_K(x):=U_{v_K}(x)=U[\mu_K](x).
\end{gather}
In such a case this potential is also equal to the probability that
the process $X_t$ starting from $x$ ever hits $K$.  

We also have 
\begin{gather}
-\lap U_K=\mu_K.  
\end{gather}

The total mass of the equilibrium measure is called the capacity 
of the set $K$ and we write 
\begin{gather}
  \Cap(K)=\int d\mu_K.
\end{gather}

The definition of the scattering length is the extension of the above to arbitrary
positive potentials. Towards this end, we first we define an analog of the capacitory measure
\begin{gather}
  \mu_v(x)=v(x)(1-U_v(x)).
\end{gather}
We want to show that $U_v(x)=U[\mu_v](x)$.
We have
\begin{gather}
  \begin{split}
    U[\mu_v](x)&
    =
    \E^x\left(\int_0^\infty v(X_s)\E^{X_s}\left( e^{-\int_0^\infty v(X_r)dr} \right)ds \right)
    \\&=
    \E^x\left(\int_0^\infty v(X_s)\E\left(e^{-\int_0^\infty v(X_r)dr}\circ \theta_s |{\mathcal F}_s \right) \right)
    \\&=
    \int_0^\infty \E^x\E\left(v(X_s)e^{-\int_0^\infty v(X_r)dr}\circ \theta_s |{\mathcal F}_s\right)ds.
  \end{split}
\end{gather}
The last equality follows from Fubini theorem and the fact that $v(X_s)$ is ${\mathcal F}_s$ measurable.
Hence 
\begin{gather}
  \begin{split}
    U[\mu_v](x)&=
    \E^x\left( \int_0^\infty v(X_s)e^{-\int_0^\infty v(X_r\circ \theta_s)dr}ds\right)
    \\&=
    \E^x\left( \int_0^\infty v(X_s)e^{-\int_s^\infty v(X_r)dr}ds\right). 
  \end{split}
\end{gather}

By Lemma \ref{finite} we have
\begin{gather}
  \int_0^\infty v(X_s)ds<\infty\;a.\;s.
\end{gather}
for almost every starting points $x\in \Rd$. 
Therefore the function
\begin{gather}
  f(s)=\int_s^\infty v(X_r)dr 
\end{gather}
is absolutely continuous for almost all paths of the process $X_s$ and
so is $e^{-f(s)}$. 
By the fundamental theorem of calculus
\begin{gather}
  \begin{split}
    U[\mu_v](x)&=\E^x\left( \int_0^\infty {d\over ds}e^{-\int_s^\infty v(X_r)dr}ds \right)
    \\&=
    \E^x\left(1-e^{-\int_0^\infty v(X_r)dr} \right)=U_v(x),
  \end{split}
\end{gather}
for almost every $x\in\Rd$. 

Note that if $v$ is bounded than by the second part of the Lemma \ref{finite} 
last equality holds for all $x$. Since $e^{-f(s)}$ is nondecreasing, its derivative
exists almost everywhere and $U_v(x)\geq U[\mu_v](x)$ for all $x\in\Rd$.

We showed that $U_v(x)$ is equal to the potential of a positive measure $\mu_v$, hence
\begin{gather}
  -\lap U_v = \mu_v.
\end{gather}
>From here we finally define the scattering length of $v$ as
\begin{equation}
  \Gamma(v)=\int_\Rd d\mu_v(x).
\end{equation}

\section{Properties of Scattering Length}

In this section we prove several useful properties of the scattering length and
the capacitory potential. First we establish some basic upper bounds for $U_v$ 
and $\Gamma(v)$.  
\begin{prop}\label{basic}
  Let $v\in L^1(\Rd)$. Then
  \begin{enumerate}
    \item $\Gamma(v)\leq \|v\|_1$,
    \item $-(\lap U_v,U_v)=\int U_vd\mu_v\leq\Gamma(v)$,
    \item if $B\subset\Rd$ bounded, then
      \begin{displaymath}
        \int_B U_v(x)dx\leq C(B)\Gamma(v).
      \end{displaymath}
  \end{enumerate}
\end{prop}
\begin{proof}
 The  first two inequalities follow from our definitions of scattering
  length and capacitory potential. For the last one we have
  \begin{displaymath}
    \begin{split}
      \int_B U_v(x)dx&=C\int_B\int_\Rd {d\mu_v(y)\over|x-y|^{d-\alpha}}dx
      \\&\leq
      C\left(\sup_y\int_B{dx\over|x-y|^{d-\alpha}}\right)\Gamma(v)
      \\&=
      C(B)\Gamma(v).
    \end{split}
  \end{displaymath}
\end{proof}
Next we prove some monotonicity and convergence properties of scattering. 
\begin{prop}\label{propconv} 
  Let $v,v_n,w\in L^1(\Rd)$ be positive. Then
  \begin{enumerate}
    \item if $v\leq w$ a.e. then $U_v\leq U_w$ a.e. and $\Gamma(v)\leq \Gamma(w)$,
    \item if $v_n(x)$ is a.e. nondecreasing and converges a.e. to $v$ then $U_{v_n}$ 
      is a.e. nondecreasing and converges a.e. to $U_v$, and
      $\Gamma(v_n)$ is nondecreasing and converges to $\Gamma(v)$.
  \end{enumerate}
\end{prop}

\begin{proof}
  Both cases for the capacitory potential follow immediately from the 
  definition for every $f\in L^1(\Rd)$.

  To prove the results about the scattering length, we need another formula 
  for $\Gamma(v)$ if $\supp v$ is bounded.  
  Consider a compact neighborhood $K$ of the support of $v$ such that 
  $v\subset\subset K$. 
  Let $U_K$ be its capacitory potential (note that $U_K=1$ on $\supp v$) and $\mu_K$ an
  equilibrium measure. 
  Then 
  \begin{gather}\label{capint}
  \Gamma(v)=\int U_K d\mu_v=\int U[\mu_K]d\mu_v=\int U[\mu_v]d\mu_K.
  \end{gather}
  Therefore monotonicity of the scattering length follows from the monotonicity of
  the potentials, if we take $K$ such that $\supp w \subset\subset K$.
  By the monotone convergence theorem part $(2)$ of the proposition is also true
 for functions of bounded support.
  Now we go back to the general case of arbitrary 
  positive $L^1(\Rd)$ functions.
  Let $v_n$  be a.e. nondecreasing with bounded supports such that $v_n$ converges a.e. to $v\in L^1(\Rd)$. 
  Suppose that $w\leq v$ a.e. and let $w_n=\min\{v_n,w\}$. We have
  \begin{gather}
    \int_\Rd v_n(x)(1-U_{v_n}(x)) dx=\Gamma(v_n)\geq\Gamma(w_n)=
\int_\Rd w_n(x)(1-U_{w_n}(x)) dx.
  \end{gather}
  Both functions under the integrals are bounded above by $v$, hence by the 
  dominated convergence theorem
  \begin{gather}
    \Gamma(v)=\int_\Rd v(x)(1-U_v(x))dx\geq \int_\Rd w(x)(1-U_w(x))dx=
    \Gamma(w).
  \end{gather}
  The second part of $(2)$ now follows from monotonicity and dominated convergence theorem.
\end{proof}

Among other interesting properties we have
\begin{prop}\label{propsum}
  For $r>0$ and $v,w\in L^1(\Rd)$ we have 
  \begin{enumerate}
    \item $U_{v+w}\leq U_v+U_w$ and $\Gamma(v+w)\leq\Gamma(v)+\Gamma(w)$ and 
    \item if $v_r(x)=r^\alpha v(rx)$, then $U_{v_r}(x)=U_v(rx)$ and
      $\Gamma(v_r)=r^{\alpha-d}\Gamma(v)$,
  \end{enumerate}
\end{prop}
\begin{proof}
The  first inequality follows from the inequality 
  $1-e^{-a-b}\leq (1-e^{-a})+(1-e^{-b})$ which is valid for any non--negative numbers $a$
  and $b$. The second inequality follows from the first one and the monotonicity of
  the potentials.

  The second part of the proposition can be easily verified by a direct calculations.
\end{proof}

Next, we consider $\Gamma(av)$, where $a$ is a positive constant. We are interested in knowing what
happens when $a\to0$ or $a\to\infty$. The following two propositions give the answer when $a\to0$. 
\begin{prop}
  Let $p^{-1}+q^{-1}=1$ and $1<p<\infty$. Assume that $v\in L^p(\Rd)$ and $\supp v\subset B$ bounded. Then
  \begin{displaymath}
    \Gamma(\eps v)=\eps\|v\|_1-O(\eps^{1+1/q}\|v\|_p\|v\|_1^{1/q}).
  \end{displaymath}
\end{prop}
\begin{proof}
  By Proposition \ref{basic} and by definition we have
  \begin{gather}
    \|U_v\|_{L^1(B)}\leq C(B)\Gamma(v),\\
    \|U_v\|_\infty\leq 1.
  \end{gather}
  Hence
  \begin{displaymath}
    \begin{split}  
    \|U_v\|_{L^q(B)}&= \left(\int_B |U_v|^q\right)^{1/q}\leq
    \left(\int_B |U_v|\right)^{1/q}\leq
    C(B)\Gamma(v)^{1/q}.
    \end{split}
  \end{displaymath}
But 
  \begin{displaymath}
    \begin{split}
      \int vU_v dx&\leq \|v\|_p\|U_v\|_{L^q(B)}
      \\&\leq
      C(B)\|v\|_p\Gamma(v)^{1/q}
      \\&\leq
      C(B)\|v\|_p\|v\|_1^{1/q}.
    \end{split}
  \end{displaymath}
  Hence,
  \begin{displaymath}
    \Gamma(\eps v)=\int \eps v(1-U_{\eps v}) dx =
   \eps\|v\|_1-O(\eps^{1+1/q}\|v\|_p\|v\|_1^{1/q}).
  \end{displaymath}
\end{proof}

\begin{prop}
  If $v\in L^1(\Rd)$, then
  \begin{displaymath}
    \lim_{\eps\searrow 0}\frac1\eps \Gamma(\eps v)=\|v\|_1.
  \end{displaymath}
\end{prop}
\begin{proof}
  By the definition of scattering length we have
  \begin{equation}
    \|v\|_1-\frac1\eps\Gamma(\eps v)=\int vU_{\eps v}dx.
  \end{equation}
  By Proposition \ref{basic}, $\|U_{\eps v}\|_{L^1(B)}\leq C(B)\|v\|_1\eps$.
  Therefore $U_{\eps v}\to0$ in measure, and the same is true for $vU_{\eps v}$. Since
  $0\leq U_{\eps v}\leq 1$, by dominated convergence theorem $\int vU_{\eps v}\to 0$ and this completes the proof. 
\end{proof}

If the potential $v$ is large, the scattering length is close to the capacity of the 
support of $v$.

\begin{prop}
Consider $v\geq0$ bounded, with compact support $K$. Let us also assume that $K$ is 
Kac regular. Then
\begin{displaymath}
  \Gamma(v)\leq\Cap K.
\end{displaymath}
If $v_i(x)\nearrow v_K(x)$, then
\begin{displaymath}
  \Gamma(v_i)\nearrow \Cap K.
\end{displaymath}
\end{prop}
\begin{proof}
  Suppose that $v\geq 0$ is supported in $K$, where $K$ is Kac regular set as described above. 
  Let $B$ be a ball such that $K$ is contained in the interior of $B$. 
  By (\ref{capint})
  \begin{displaymath}
    \Gamma(v)=\int U_v(x)d\mu_B(x),
  \end{displaymath}
  and
  \begin{displaymath}
    \Cap K=\int U_K(x)d\mu_B(x).
  \end{displaymath}
  But $U_v\leq U_K$, so $\Gamma(v)\leq \Cap K$.
 The  second part of the proposition follows from monotone convergence theorem and
  from the first part.
\end{proof}

\section{Upper bound of Theorem \ref{main}}
 
 In this section we prove the upper bound in the main result using variational
 characterization of eigenvalues. First eigenvalue $\lambda_1(v)$ can be calculated using
 Rayleigh quotient
 \begin{gather}\label{ratio}
   \lambda_1(v)=\inf_{\varphi\in W^{\alpha/2,2}(\Omega)}
   \frac{\mathcal{E}_Y(\varphi,\varphi)+\int_\Omega v\varphi^2}{\int_\Omega \varphi^2}
 \end{gather}
 Our strategy is  to choose a function $\varphi$ which will give desired bound.
Our claim is that the function  $\varphi=U_v-1$ will do the job. This function is in $W^{\alpha/2,2}(\Omega)$, since $U_v$
 is in the domain of $\lapN$ and $\Omega$ is bounded.  We have 
 \begin{displaymath}
   \begin{split}
     \mathcal{E}_Y(\varphi,\varphi)&+\int_\Omega v\varphi^2=
     \mathcal{E}_Y(U_v,U_v)+\int_\Omega v\varphi^2\\
     &=
     \frac12\int_{\Omega\times\Omega}\frac{(U_v(x)-U_v(y))^2}{|x-y|^{d+\alpha}}dxdy
     +\int_\Omega v\varphi^2\\
     &\leq
     \frac12\int_{\Rd\times\Rd}\frac{(U_v(x)-U_v(y))^2}{|x-y|^{d+\alpha}}dxdy
     +\int_\Omega v\varphi^2\\
     &=
     \mathcal{E}_X(U_v,U_v)+\int_\Omega v\varphi^2=
     \int_\Rd (-U_v\Delta^{\alpha/2}U_v+v(U_v-1)^2)\\
     &=
     \int_\Rd (vU_v(1-U_v)+vU_v(U_v-1)-v(U_v-1))=\Gamma(v).
   \end{split}
 \end{displaymath}
 On the other hand,
 \begin{gather}
   \begin{split}
     \int_\Omega \varphi^2\geq |\Omega|-2C(\Omega)\Gamma(v).
   \end{split}
 \end{gather}
 Hence if $\Gamma(v)\leq |\Omega|/(4C(\Omega))=\beta(\Omega)$ the last expression is 
 comparable to the volume of $\Omega$.
 This completes the proof of the upper bound.

 \section{Comparison lemma}
In this section we define a class of processes called stable-like processes. Then
we pick a suitable representative to be used in the proof of the lower bound in 
the main result. Two key lemmas used in the next section are also given.
 
 The processes $X_t$ and $Y_t$ are examples of a larger class of processes defined in
 \cite{CK}, called stable-like processes $Z_t$.  These
 processes have generators with quadratic forms 
 \begin{equation}
   {\mathcal E}_Z(u,u)=\int_\Omega\int_\Omega c(x,y){(u(x)-u(y))^2\over|x-y|^{d+\alpha}}dxdy.
 \end{equation}
 Here  $c(x,y)$ is a symmetric function satisfying $0<c\leq c(x,y)\leq C<\infty$ for all $x$, $y$ where $c$ and $C$ are constants independent of $x$ and $y$. 
The domain of this form is the same as the domain
 of the ``Neumann'' fractional Laplacian, 
 namely $W^{\alpha/2,2}(\Omega)$. For more detail about this class we 
 refer the reader to \cite{CK}. 

 Let $B_t$ be a Brownian motion running at twice the usual speed, and 
 $U_t$ be a reflected 
 Brownian motion in a cube.  That is,  $U_t$ is the process generated by 
 the Laplacian with  Neumann boundary conditions in the cube. 
 We will use a subordination technique (see \cite{Sa}) to obtain stable 
 processes from these processes.
 Let  $A_t$ be a positive $\alpha/2$-stable subordinator  independent
 of $B_t$ and $U_t$.
 If we subordinate a Brownian motion with $A_t$ we get a $\alpha$-stable
 process. In other words $X_t=B_{A_t}$. Let $V_t$ be
 the process $U_t$ subordinated with the same subordinator $A_t$. The resulting 
 process is a stable-like process (see \cite{CK}).  However,  it is not the same as a the reflected
 stable process $Y_t$ (see \cite{BBC}).

The following lemma gives a comparison between expected values of the multiplicative potentials
of $X_t$ and $V_t$
 \begin{lem}\label{proccomp}
   Let $\supp(v)\subset\Omega$, where $\Omega$ is a cube.  Then 
   \begin{equation}
     \E^x\left\{ \exp\left( -\int_0^t v(V_s)ds \right) \right\}
     \leq 
     \E^x\left\{ \exp\left( -\int_0^t v(X_s)ds \right) \right\}
   \end{equation}
 \end{lem}
 \begin{proof}
Define $g$ as follows
    \begin{equation}
      g(x)=
      \begin{cases}
        x-2n\mbox{ ,if } x\in[2n,2n+1),\;n\in Z,\\
        2n-x\mbox{ ,if } x\in[2n-1,2n),\; n\in Z.
      \end{cases}
    \end{equation}
One can think about $g$ as a function that continuously folds a real line
into a unit interval.
    Using a tensor product we can define  
    $$f(x_1,x_2,\cdots,x_d)=g(x_1)\otimes g(x_2)\otimes \cdots
    \otimes g(x_d).$$ 
    We have $U_t=f(B_t)$ for a reflected Brownian motion $U_t$ on $[0,1]^d$.
    Since $1$-dimensional components of $B_t$ (and $U_t$ on a cube) are independent of
    each other and are transition invariant we have
    \begin{gather}
      U_t=f(B_t),
    \end{gather}
    for arbitrary cube.

    This gives:
    \begin{displaymath}
      V_t=U_{A_t}=f(B_{A_t})=f(X_t).
    \end{displaymath}
    Now we can define $\tilde{v}(x)=v(f(x))$ so that $\tilde{v}=v$ on 
    $\supp(v)$.
    We have:
    \begin{displaymath}
      \begin{split}
     \E^x&\left\{ \exp\left( -\int_0^t v(V_s)ds \right) \right\}=
     \E^x\left\{ \exp\left( -\int_0^t v(f(X_s))ds \right) \right\}\\
     &= 
     \E^x\left\{ \exp\left( -\int_0^t \tilde{v}(X_s)ds \right) \right\}
     \leq 
     \E^x\left\{ \exp\left( -\int_0^t v(X_s)ds \right) \right\}
      \end{split}
    \end{displaymath}
 \end{proof}
 The first eigenvalues of the Schr\"odinger operators of the generators 
 of two arbitrary stable-like processes are comparable. In particular we have 
 \begin{lem}\label{eigcomp}
   Let $\lambda_1^V$ be the first eigenvalue of the operator $-A+v$, where 
   $A$ is the generator of $V_t$. Let also $\lambda_1(v)$ be as in Theorem \ref{main}, i.e. the
   first eigenvalue of the Schr\"odinger operator for the ``Neumann'' fractional 
   Laplacian. Then
   \begin{gather}
     c\lambda_1^V\leq \lambda_1(v)\leq C\lambda_1^V,
   \end{gather}
   where $c$ and $C$ are positive constants.
 \end{lem}

 \begin{proof}
   Process $V_t$ is a stable-like process, hence
   \begin{gather}
   {\mathcal E}_V(u,u)=\int_\Omega\int_\Omega c_V(x,y){(u(x)-u(y))^2\over|x-y|^{d+\alpha}}dxdy.
   \end{gather}
   and $c\leq c_V(x,y)\leq C$ for some constants $c$ and $C$. We can assume that
   $c<1$ and $C>1$.
   Given any positive potential $v$
   \begin{gather}
       \frac1C\left({\mathcal E}_V(u,u)+\int vu^2\right)\leq 
       {\mathcal E}_Y(u,u)+\int vu^2\leq
       \frac1c\left({\mathcal E}_V(u,u)+\int vu^2\right). 
   \end{gather}
   By (\ref{ratio}) we get the inequality between the eigenvalues.
 \end{proof}

 \section{Lower bound of Theorem \ref{main}}
 Consider process $V_t$ defined in \S5. By  Lemma \ref{eigcomp}, it is 
 enough to prove the lower bound for this process.
Let $A$ be its generator. It is enough to prove that there exists $t$ such that 
\begin{gather}
||e^{t(A-v)}||_2\leq e^{-C\Gamma(v)}.
\end{gather}

We will prove this inequality using heat kernels associated with various processes.
Let $u_A(t,x,y)$ be a heat kernel associated with the Schr\"odinger operator of
the operator $A$, i.e. a function satisfying
\begin{gather}
  e^{t(A-v)}f(x)=\int_\Omega u_A(t,x,y)f(y)dy,
\end{gather}
for every bounded $f$.
We have to show that 
\begin{displaymath}
  \int_\Omega u_A(x,y,t)dy\leq e^{-C\Gamma(v)},
\end{displaymath} 
Using Feynman-Kac formula we get
\begin{gather}
  \int_\Omega u_A(t,x,y)dy=\E^x\left(e^{-\int_0^\infty v(V_s)ds}\right).
\end{gather}
Let $u(t,x,y)$ be a heat kernel associated with the Schr\"odinger operator
of the fractional Laplacian. Then
\begin{gather}
  \int_\Rd u(t,x,y)dy=\E^x\left(e^{-\int_0^\infty v(X_s)ds}\right).
\end{gather}
By the Lemma \ref{proccomp} it is now enough to prove that
\begin{displaymath}
  1-U_v^t(x)=\int_\Rd u(t,x,y) dy\leq e^{-C\Gamma(v)}.
\end{displaymath}
First we need an upper bound for the capacitory potential
\begin{gather}\label{Ubound}
  U_v(x)\geq C\int \frac{d\mu_v(y)}{|x-y|^{d-\alpha}}
  \geq C\Gamma(v)(\diam(\Omega))^{\alpha-d}.
\end{gather}
Using the semigroup property of $u(t,x,y)$
\begin{gather*}
  \begin{split}
    \int u(t,x,y)&U^s_v(y)dy=\int u(t,x,y)\left(1-\int u(y,z,s)dz\right)dy\\
    &=
    \int u(t,x,y) dy-\int u(t+s,x,z)dz=U^{t+s}_v(x)-U^t_v(x),
  \end{split}
\end{gather*}
If we let $s$ tend to $\infty$ we get
\begin{gather}\label{difference} 
    U_v(x)-U^t_v(x)=\int u(t,x,y) U_v(y) dy.
\end{gather}
Let $p(t,x,y)$ be a heat kernel associated with the process $X_t$.
Since our potentials $v$ are nonnegative, we have
\begin{gather}
  u(t,x,y)\leq p(t,x,y).
\end{gather}
Using this inequality we obtain
\begin{displaymath}
  \begin{split}
    U_v(x)-U^t_v(x)&=\int u(t,x,y) U_v(y) dy\leq
    \int p(t,x,y)U_v(y)dy\\
    &=C\int\int \frac{d\mu_v(z)}{|z-y|^{d-\alpha}}p(t,x,y)dy\\
    &=C\int\int \frac{p(t,x,y)}{|z-y|^{d-\alpha}}dy d\mu_v(z)\\
    &\leq
    C\Gamma(v)\sup_{z\in\Omega,\;x\in\Rd}\int \frac{p(t,x,y)}{|z-y|^{d-\alpha}}dy,
  \end{split}
\end{displaymath}
We need to show that supremum tends to 0 as $t$ tends to $\infty$. 
Then we can take $t_0$ large enough so that
\begin{gather}
  U_v-U_v^{t_0}\leq U_v/2.
\end{gather}

And using (\ref{Ubound})
\begin{displaymath}
  1-U^{t_0}_v(x)\leq 1-U_v(x)/2\leq 1-C/2\Gamma(v)(\diam\Omega)^{\alpha-d}=
  e^{-C\Gamma(v)}.
\end{displaymath}

The only thing left to prove is the following
\begin{lem} Let $p(t,x,y)$ be a heat kernel associated with the process  $X_t$. 
  Then
  \begin{displaymath}
    \lim_{t\To\infty}\sup_{z,x\in\Rd} 
    \int_\Rd \frac{p(t,x,y)}{|y-z|^{d-\alpha}}dy = 0.
  \end{displaymath}
\end{lem}
\begin{proof}
  We divide $\Rd$ into $3$ sets, and estimate the integral on each of these sets.

Let $B_1=\left\{ y\in\Rd:\;|y-z|\leq t^{1/\alpha} \right\}$. Since 
  $p(t,x,y)\leq c t^{-d/\alpha}$, we have
  \begin{displaymath}
    \begin{split}
      \int_{B_1} \frac{p(t,x,y)}{|y-z|^{d-\alpha}}dy&\leq
      \int_{B_1} \frac{c t^{-d/\alpha}}{|y-z|^{d-\alpha}}dy\\
      &=
      c t^{-d/\alpha}\int_0^{t^{1/\alpha}} \frac{r^{d-1}}{r^{d-\alpha}}dr\\
      &=
      c t^{-d/\alpha} t=c t^{1-d/\alpha}.
    \end{split}
  \end{displaymath}
  But $d>\alpha$, so the last expression tends to 0 when $t$ tends to $\infty$.

Let  
  $B_2=\left\{ y\in\Rd:\;|y-z|>t^{1/\alpha},\; |x-y|\leq t^{1/\alpha} \right\}$.
  By the same estimate for $p(t,x,y)$
  \begin{displaymath}
    \begin{split} 
      \int_{B_2} \frac{p(t,x,y)}{|y-z|^{d-\alpha}}dy&\leq
      \int_{B_2} \frac{c t^{-d/\alpha}}{t^{(d-\alpha)/\alpha}}dy\\
      &=
      c|B_2| t^{-d/\alpha-(d-\alpha)/\alpha}=
      c t^{-(d-\alpha)/\alpha}
    \end{split}
  \end{displaymath}
  Hence this part of the integral also tends to 0.

Let 
  $B_3=\left\{ y\in\Rd:\;|y-z|>t^{1/\alpha},\; |x-y|> t^{1/\alpha} \right\}$. Here
  we can use the fact that $p(t,x,y)\leq c \frac{t}{|x-y|^{d+\alpha}}$, to get
  \begin{displaymath}
    \begin{split}
      \int_{B_3} \frac{p(t,x,y)}{|y-z|^{d-\alpha}}dy&\leq
      \int_{B_3} \frac{c t}{|x-y|^{d+\alpha}t^{(d-\alpha)/\alpha}}dy\\
      &=
      c t^{1-(d-\alpha)/\alpha}\int_{t^{1/\alpha}}^\infty 
      \frac{r^{d-1}}{r^{d+\alpha}}dr\\
      &=
      c t^{1-(d-\alpha)/\alpha}(t^{1/\alpha})^{-\alpha}=
      c t^{-(d-\alpha)/\alpha}.
    \end{split}
  \end{displaymath}
  Therefore each part of the integral tends to 0 as $t$ tends to $\infty$. This completes  the proof of
Theorem \ref{main}.
\end{proof}

\section*{Acknowledgements}
The Author would like to thank Professor Rodrigo Ba\~nuelos, his Ph. D. thesis advisor, for
guidance on this paper. The author would also like to thank Professor Krzysztof Bogdan
for helpful discussions about the definition of the scattering length.

\bibliographystyle{amsplain}

\end{document}